\documentclass{article}
\title{Torsion in almost K\"ahler geometry 
}
\author{Paul-Andi Nagy}
\date{\today}
\usepackage{amsfonts,amssymb}
\usepackage{minitoc_href}

\newtheorem{teo}{Theorem}[section]
\newtheorem{lema}{Lemma}[section]
\newtheorem{pro}{Proposition}[section]
\newtheorem{defi}{Definition}[section]
\newtheorem{rema}{Remark}[section]
\newtheorem{coro}{Corollary}[section]
\newtheorem{nr}{}[section]
\begin{document}
\maketitle
\abstract{\normalsize We study almost K\"ahler manifolds whose curvature tensor satisfies the second 
curvature condition of Gray (shortly ${\cal{AK}}_2$). This condition is interpreted in terms of the first canonical 
Hermitian connection. 
It turns out that this condition forces the torsion of this connection to be parallel in directions orthogonal to the K\"ahler nullity of the almost 
complex structure. We prove a local structure result for ${\cal{AK}}_2$ manifolds, showing that the basic pieces are manifolds 
with parallel torsion and special almost K\"ahler manifolds, a class generalizing, to some algebraic extent, the class of $4$-dimensional 
${\cal{AK}}_2$-manifolds. In the case of parallel torsion, the Einstein condition and the reducibility of 
the canonical Hermitian connection is studied.
\large
\tableofcontents
\section{Introduction}
An  almost K\"ahler manifold (shortly ${\cal{AK}}$) is a Riemannian manifold $(M^{2n},g)$, together with a compatible almost complex structure 
$J$, such that the K\"ahler form $\omega=g(J \cdot, \cdot)$ is closed. Hence, almost K\"ahler geometry is nothing else that symplectic 
geometry with a prefered metric and complex structure. Since symplectic manifolds often arise in this way is rather 
natural to ask under which conditions on the metric we get integrability of the almost complex $J$. In this direction, 
a famous conjecture of S. I. Golberg asserts that every compact, Einstein, almost K\"ahler manifold is, in fact, K\"ahler. At our present 
knowledge, this conjecture is still open. Nevertheless, they are a certain number of partial results, supporting this 
conjecture. First of all, K. Sekigawa proved \cite{Seki1} that 
the Goldberg conjecture is true when the scalar curvature is positive.
We have to note that the Golberg conjecture is definitively not true with the compacity assumption removed. In fact, there are 
Hermitian symmetric spaces of non-compact type of any complex dimension $n \ge 3$ admitting almost K\"ahler structures commuting with the invariant 
K\"ahler one \cite{Apo2}. Not that, at the opposite, the real hyperbolic space of dimension at least $4$ do not admits, even locally, orthogonal almost K\"ahler structures 
\cite{Olszak, arm2}. In dimension $4$, examples of local Ricci flat almost K\"ahler metrics are constructed in 
\cite{Apo, arm2,Nur} . 
In the same paper, a potential source of compact almost K\"ahler, Einstein manifolds is considered, 
namely those compact K\"ahler manifolds whose Ricci tensor admits two distinct, constant eigenvalues; integrability 
is proven under certain positivity conditions. 
The rest of known results, most of them 
enforcing or replacing the Einstein condition with some other natural curvature assumption are mainly 
in dimension $4$. To cite only a few of them, we mention the beautiful series of papers \cite{Apo1, Apo3, Apo4} giving 
a complete local and global classfication of almost K\"ahler manifolds of $4$ dimensions satisfying the second and third 
Gray condition on the Riemannian curvature tensor. Other recent results, again in $4$-dimensions, are concerned with the study of local obstructions 
to the existence of Einstein metrics \cite{arm2}, $\star$-Einstein metrics \cite{Seki2}, etc.
\par
In this paper our main object of study will be the class of almost K\"ahler manifolds satisfying the second curvature 
condition of Gray. To the best of our knowledge, the only classification available at present is that of \cite{Apo1} where 
it is shown that every $4$-dimensional almost K\"ahler manifold satisfying the second curvature condition of Gray is 
locally isometric to the unique $3$-symmetric space of $4$ dimensions. Some classes of examples are also 
known, such as twistor spaces over quaternionic-K\"ahler manifolds of negative scalar curvature 
\cite{Davidov, Ivanov}. \par
Our approach to the study of the class of almost K\"ahler manifolds satisfying the second Gray curvature condition 
(shortly ${\cal{AK}}_2$) will be directed from the point of view of the canonical Hermitian 
connection. Actually, we are going to study the geometric as well as the algebraic effects of the second curvature 
condition of Gray over the torsion of the last mentioned connection. 
The main result of the present paper, giving a local structure theorem concerning 
almost K\"ahler manifolds in the class ${\cal{AK}}_2$ is the following 
\begin{teo}
Let $(M^{2n},g,J)$ be an almost K\"ahler manifold in the class ${\cal{AK}}_2$. Let $U$ be an open set where the 
K\"ahler nullity has constant rank. Then there is an open dense (with respect to the induced topology) set $D$ in $U$ such that 
around each point of $D$ the manifold  $M$ is locally the Riemannian product of a almost K\"ahler manifold whose 
first canonical connection has parallel torsion and a special ${\cal{AK}}_2$-manifold.
\end{teo}
The precise definition of special ${\cal{AK}}_2$ manifolds is given at the end of the section 4. They are those 
supporting almost K\"ahler structures for which the integral manifolds of the 
distribution orthogonal to the K\"ahler nullity are K\"ahler, with respect to the induced structure. Note that 
algebraically (see definition 4.1 for details) this property is automatically satisfied in $4$-dimensions. \par
Theorem 1.1 shows that the study of the torsion of an almost K\"ahler manifold of class ${\cal{AK}}_2$ reduces, in the local 
sense precised below, to the study of the structure of the torsion of a special almost K\"ahler manifold. \\ \par 
Concerning almost K\"ahler manifolds with parallel torsion, and in connection with the existence problem of Einstein almost-K\"ahler metrics, we are able to prove the following: 
\begin{teo}
For every almost K\"ahler manifold with parallel torsion the holonomy representation 
of the canonical Hermitian connection is reducible, in the real 
sense. Furthermore, if such a manifold is Einstein, then it has to be K\"ahler. 
\end{teo}
It follows that an Einstein ${\cal{AK}}_2$ manifold is locally the product of a K\"ahler Einstein manifold and an Einstein special 
${\cal{AK}}_2$-manifold. Note the difference with nearly-K\"ahler manifolds where many Einstein homogenenous examples exist (see \cite{Nagy2}).\par
Our paper is organised as follows. In section 2 we recall some general, well known facts of almost K\"ahler 
geometry. In section 3, the main technical ingredient of this paper is proved: using the first canonical Hermitian 
connection we give an interpretation of the second Gray condition on curvature in terms of the torsion of 
the last mentioned connection. Namely, we show that the associated (bundle valued) $1$-form has to be closed. 
We study this condition using some standard methods and we prove that the torsion of the canonical Hermitian 
connection has to be parallel in directions orthogonal to the K\"ahler nullity of the almost complex 
structure. Note that this result continues to hold in the more general case of quasi-K\"ahler manifolds 
satisfying the second Gray condition on curvature. \par
In section $4$ we obtain a preliminary decomposition result of the distribution orthogonal to the K\"ahler nullity. 
We use the first Bianchi identity for the canonical Hermitian connection coupled with the almost K\"ahler 
condition and in presence 
of the partial parallelism established in section $3$ in order to obtain informations about the algebraic structure of the 
torsion. Note that in the case of nearly K\"ahler geometry, and for an arbitrary holonomy decomposition this approach 
was sufficient to extract all necessary algebraic information about the torsion of the canonical Hermitian connection 
(see \cite{Nagy2}). The additional difficulties one has to face in the ${\cal{AK}}_2$ case are due on the one hand to the fact 
that we have only a partial parallelism for the torsion and on the other hand to the fact that the almost K\"ahler  
condition (viewed algebraically on the first jet of the almost complex structure) is more complicated than the nearly 
K\"ahler one. \par
The decomposition mentioned below is 
brought to its final form, leading to the proof of theorem 1.1 in section $5$. Here we start from the second 
Bianchi identity for 
the canonical Hermitian connection in order to obtain structural properties of the Hermitian 
curvature tensor. The latter are used to 
decide which parts of the decomposition of the distribution orthogonal to the K\"ahler nullity 
arising from the the algebraic study of section $4$ cannot geometrically occur. \par
Finally, in section 6, a proof of the theorem 1.2 is given, by applying Sekigawas's formula (in fact its pointwise 
version develloped in \cite{Apo2}) to the particular case of parallel torsion. \par
\section{Preliminaries}
Let us consider an almost Hermitian manifold $(M^{2n},g,J)$, that is a Riemannian manifold 
endowed with a compatible almost complex structure. We denote by $\nabla$ the Levi-Civita 
connection of the Riemannian metric $g$. Consider now the tensor $\nabla J$, the first derivative 
of the almost complex structure, and recall that for all $X$ in $TM$ we have that 
$\nabla_XJ$ is a skew-symmetric (with 
respect to $g$) endomorphism of $TM$, which anticommutes with $J$. The tensor 
$\nabla J$  can be used to distinguish various classes of almost 
Hermitian manifolds. For example, $(M^{2n},g,J)$ is quasi-K\"ahler iff 
$$\nabla_{JX}J=-J\nabla_XJ $$
for all $X$ in $TM$. If $\omega=g(J \cdot, \cdot)$ denotes the K\"ahler form of the 
almost Hermitian structure $(g,J)$, we have an almost K\"ahler structure 
iff $d \omega=0$. We also recall the well known fact that almost K\"ahler manifolds 
are always quasi-K\"ahler. \par
The almost complex structure $J$ defines a Hermitian structure if it is integrable, that is 
the Nijenhuis tensor $N_J$ defined by 
$$ N_J(X,Y)=[JX,JY]-[X,Y]-J[X,JY]-J[JX,Y]$$
for all vector fields $X$ and $Y$ on $M$ vanishes. This is also equivalent to 
$$ \nabla_{JX}J=J\nabla_XJ$$
whenever $X$ is in $TM$. Therefore, an almost K\"ahler manifold which is also 
Hermitian must be K\"ahler. \par
In this paper we will deal mainly with almost K\"ahler (${\cal{AK}}$ for short)-manifolds, altough we will autorize us short 
excursions to the quasi-K\"ahler class. We begin to recall some basic facts about the various notions 
of Ricci tensors. In the rest of this section $(M^{2n},g,J)$ will be an ${\cal{AK}}$ manifold. \par
Let $Ric$ be the Ricci tensor of the Riemannian metric $g$. We denote by $Ric^{\prime}$ and 
$Ric^{\prime \prime}$ the $J$-invariant resp. the $J$-anti-invariant part of the tensor $Ric$. Then the 
Ricci form is defined by 
$$ \rho=<Ric^{\prime}J \cdot, \cdot>.$$
We define the $\star$-Ricci form by 
$$ \rho^{\star}=\frac{1}{2} \sum \limits_{i=1}^{2n}R(e_i, Je_i)$$
where $\{ e_i, 1 \le i \le 2n \}$ is any local orthonormal basis in $TM$. Note that $\rho^{\star}$ is not, in general,
$J$-invariant. The $\star$-Ricci form is related 
to the Ricci form by 
\begin{nr} \hfill 
$ \rho^{\star}-\rho=\frac{1}{2} \nabla^{\star} \nabla \omega. \hfill $
\end{nr}
The (classical) proof of this fact consists in using the Weitzenb\"ock formula for the harmonic 
$2$-form $\omega$. Taking the scalar product with $\omega$ we obtain : 
$$ s^{\star}-s=\frac{1}{2}\vert \nabla J\vert^2$$
where the $\star$-scalar curvature is defined by $s^{\star}=2<R(\omega), \omega>$
\par 
We now come to properties of the curvature tensor of an almost K\"ahler manifold. Recall 
first that we have the decomposition : 
$$ \Lambda^2(M)=\Lambda^{1,1}(M) \oplus \Lambda_{-}^2(M)$$
into $J$-invariant and $J$-anti-invariant parts, where the action of $J$ on a two form $\alpha$ is given by 
$(J\alpha)(X,Y)=\alpha(JX,JY)$ for all $X,Y$. Then the real vector bundle 
$\Lambda_{-}^2(M)$ has a complex structure ${\cal{J}}$ defined by $({\cal{J}}\alpha)(X,Y)=
-\alpha(JX,Y)$. \par
Let $R$ be the curvature tensor of the metric $g$, with the convention that $R(X,Y)=\nabla_{[X,Y]}-
[\nabla_X,\nabla_Y]$ for all vector fields $X$ and $Y$ on $M$. Let now $\tilde{R}$ be the component 
of $R$ acting trivially on $\Lambda^{1,1}(M)$. It is defined by : 
$$ \tilde{R}(X,Y,Z,U)=\frac{1}{4}(R(X,Y,Z,U)-R(JX,JY,Z,U)-R(X,Y,JZ,JU)+
R(JX,JY,JZ,JU)). $$
Now, $\tilde{R}$ decomposes further as $\tilde{R}=\tilde{R}^{\prime}+\tilde{R}^{\prime \prime}$ where 
$\tilde{R}^{\prime}$ and $\tilde{R}^{\prime \prime}$ are the components of $\tilde{R}$ commuting, resp. 
anticommuting with ${\cal{J}}$. Explicitely, we have : 
$$\begin{array}{rr}
\tilde{R}^{\prime \prime}(X,Y,Z,U)=\frac{1}{8} \biggl ( R(X,Y,Z,U)-R(JX,JY,Z,U)-R(X,Y,JZ,JU)+
R(JX,JY,JZ,JU)\\
-R(X,JY,Z,JU)-R(JX,Y,Z,JU)-R(X,JY,JZ,U)-R(JX,Y,JZ,U)
\end{array} $$ 
and 
$$\begin{array}{rr}
\tilde{R}^{\prime }(X,Y,Z,U)=\frac{1}{8} \biggl ( R(X,Y,Z,U)-R(JX,JY,Z,U)-R(X,Y,JZ,JU)+
R(JX,JY,JZ,JU)\\
+R(X,JY,Z,JU)+R(JX,Y,Z,JU)+R(X,JY,JZ,U)+R(JX,Y,JZ,U)
\biggr )
\end{array} $$ 
for all $X,Y,Z,U$ in $TM$. An important information about this tensors is given by an identity of A. Gray from  
\cite{Gray1}, page 604, Cor.4.3 : \\
\begin{nr} \hfill
$ \tilde{R}^{\prime }(X,Y,Z,U)=-\frac{1}{4}<(\nabla_XJ)Y-(\nabla_YJ)X, 
(\nabla_ZJ)U-(\nabla_UJ)Z>.
\hfill $
\end{nr}
We will now give an important formula, which can be interpreted as an obstruction to the 
existence of almost K\"ahler, non-K\"ahler structures. 
\begin{pro}\cite{Apo2}
Let $(M^{2n},g,J)$ be an almost K\"ahler manifold. Then the following holds :  
$$ \begin{array}{cc}
\Delta(s^{\star}-s)=&-4\delta (J\delta(JRic^{\prime \prime}))+8\delta(<\rho^{\star}, 
\nabla_{\cdot}\omega>)+2\vert Ric^{\prime \prime}\vert^2  \vspace{2mm}\\
& \! \! \! \! \! \! \! \! \! \! \! \! \! \! \! \!+4<\rho, \Phi-\nabla^{\star}\nabla \omega>-
\vert \Phi \vert^2-\vert \nabla^{\star}\nabla \omega\vert^2-8\vert \tilde{R}^{\prime \prime}\vert^2. \\
%&  \! \! \! \! \! \! \! \! \! \! \! \! \! \! \! \! \! \! \! \! \! \! \! \! \!-8\vert \tilde{R}^{\prime}\vert^2
%+4 \sum \limits_{i=1}^{2m}<
%\tilde{R}^{\prime}, \nabla_{e_i}\omega \otimes \nabla_{e_i}\omega>
\end{array} $$
Here, the semi-positive $2$-form $\Phi$ is defined by 
$$ \Phi(X,Y)=<\nabla_{JX}\omega, \nabla_Y \omega>.$$
and $\delta$ denotes co-differentiation with respect to $\nabla$, acting on $1$-forms and $2$-tensors.
\end{pro}
\section{Gray's curvature conditions}
This section is dedicated to interpret some of the well known conditions on the curvature tensor of an almost Hermitian 
manifold in terms of the torsion of the first canonical Hermitian connection. We begin by recalling how one 
can distinguish several classes of almost Hermitian manifolds by "the degree of ressemblance" of their 
Riemannian curvature tensor with the curvature tensor of a K\"ahler manifold. \par
Let $(M^{2n},g,J)$ be an almost Hermitian manifold. Let $R$ be the Riemannian curvature tensor 
of the metric $g$. Then the following classes of almost Hermitian manifolds appear in a natural 
way \cite{Gray1}: \\
$\begin{array}{lr} 
(G_1) : \ R(X,Y,JZ,JU)=R(X,Y,Z,U) \\
(G_2) : \  R(X,Y,Z,U)-R(JX,JY,Z,U)=R(JX,Y,JZ,U)+R(JX,Y,Z,JU) \\
(G_3) : \ R(JX,JY,JZ,JU)=R(X,Y,Z,U)
\end{array} $
$\\$
Using the first Bianchi identity it is a simple exercise to see that $G_1 \Rightarrow G_2 \Rightarrow G_3$. It is also clear that a K\"ahler 
structure satisfies all the three conditions.  Let us set now some notations. \par
Following \cite{Gray1}, let ${\cal{AK}}$ be the class of almost K\"ahler manifolds. Then the class ${\cal{AK}}_i, 1 \le i \le 3$ contains those 
almost K\"ahler manifolds whose curvature tensor satisfies the condition $(G_i)$. Obviously, we have the inclusions : 
$$ {\cal{AK}}_1 \subseteq {\cal{AK}}_2 \subseteq {\cal{AK}}_3. $$
Note that it was shown in \cite{Golberg} that locally $ {\cal{AK}}_1={\cal{K}}$, where ${\cal{K}}$ denotes the class of K\"ahler manifolds. The other 
inclusions are strict in dimensions $2n \ge 6$, as shows the examples of 
\cite{Davidov}, multiplied by K\"ahler manifolds. In the same spirit, the class ${\cal{AH}}_i, 1 \le i \le 3$ contains 
those almost Hermitian manifolds whose Riemannian curvature tensor satisfies condition $(G_i)$.  
$\\$
\par
Let us consider the first canonical connection of the almost Hermitian manifold $(M^{2n}, g, J)$ to be defined by : 
$$ \overline{\nabla}_XY=\nabla_XY+\eta_XY$$
whenever $X,Y$ are vector fields on $M$, where $\nabla$ is the Levi-Civita connection of $g$ and where, 
to save space, we setted $\eta_XY=\frac{1}{2}(\nabla_XJ)JY$. We obtain a metric Hermitian connection on 
$M$, that is $\overline{\nabla}g=0$ and $\overline{\nabla}J=0$. Recall that for each $X$ in $TM$,
$\eta_X$ is a skew-symmetric, $J$-anticommuting endomorphism of $TM$. \par
The torsion tensor of the canonical Hermitian canonical connection, to be denoted by $T$ is given by 
$$T_XY=\eta_XY-\eta_YX$$
for all $X,Y$ in $TM$. Then by the {\bf{torsion}} of the almost Hermitian manifold $(M^{2n},g,J)$ we will mean 
simply the torsion tensor of the canonical Hermitian connection. \par
For the almost Hermitian $(M^{2n},g,J)$ be almost K\"ahler one requires that the K\"ahler form 
$\omega(X,Y)=<JX, Y>$ be closed. In our notations, this is equivalent to have 
\begin{nr} \hfill
$<T_XY,Z>=-<\eta_ZX,Y> \hfill $
\end{nr}
for all $X, Y$ and $Z$ in $TM$. 
In the almost K\"ahler context, this relation will be used almost implicitely in the rest of this paper. \par  
Let $\overline{R}(X,Y)=-[\overline{\nabla}_X, 
\overline{\nabla}_Y]+\overline{\nabla}_{[X,Y]}$ be the curvature tensor of the connection 
$\overline{\nabla}$. Now a standard calculation involving the definitions yields to 
\begin{nr} \hfill 
$ \overline{R}(X,Y)Z=R(X,Y)Z+[\eta_X,\eta_Y]Z-\biggr [ d_{\overline{\nabla}}u(X,Y) \biggl ] Z \hfill $
\end{nr}
where 
$$ \biggl [ d_{\overline{\nabla}}u(X,Y) \biggr ] Z =(\overline{\nabla}_X\eta)(Y,Z)-(\overline{\nabla}_Y\eta)(X,Z)+
\eta_{T_XY}Z$$
for all vector fields $X,Y,Z$ on $M$. Note that $d_{\overline{\nabla}}u(X,Y)$ is a $J$-anticommuting endomorphism 
of $TM$, whenever $X,Y$ are tangent vectors to $M$.
\begin{rema}
In the formula (3.2), the notation $u$ stands for the tensor $\eta$, considered as a $1$-form with values in the 
bundle $\Omega^2(M)$. Then $d_{\overline{\nabla}}$, with its expression given below, is the twisted differential 
acting on twisted one forms, when considering the tangent bundle of $M$ endowed with the connection 
$\overline{\nabla}$. Since our discussion is intended to be self-contained and at the elementary level, we will 
keep things at the level of the notation after (3.2).
\end{rema}
The fact that $\overline{\nabla}$ is a Hermitian connection implies that $\overline{R}(X,Y,JZ,JU)=\overline{R}(X,Y,Z,U)$. Using this in formula (3.2), together with 
the skew-symmetry of the $J$-anticommuting endomorphism $\eta_X$ gives us : 
\begin{nr} \hfill 
$ \begin{array}{lr}
R(X,Y,Z,U)-R(X,Y,JZ,JU)=2\biggl [ (d_{\overline{\nabla}}u)(X,Y) \biggr ](Z,U).
\end{array}
\hfill $
\end{nr}
By the symmetry property of the Riemannian curvature we also deduce that 
\begin{nr} \hfill 
$ \begin{array}{lr}
R(X,Y,Z,U)-R(JX,JY,Z, U)=2\biggl [ (d_{\overline{\nabla}}u)(Z,U) \biggr ](X,Y).
\end{array}
\hfill $
\end{nr}
\begin{lema}The almost Hermitian manifold $(M^{2n},g,J)$ satisfies the condition $(G_3)$ iff 
\begin{nr} \hfill 
$ \biggl [ (d_{\overline{\nabla}}u)(Z,U) \biggr ](X,Y)=
\biggl [ (d_{\overline{\nabla}}u)(X,Y) \biggr ](Z,U).\hfill $
\end{nr}
Moreover, we have $(d_{\overline{\nabla}}u)(JX,JY)+(d_{\overline{\nabla}}u)(X,Y)=0$, hence 
$d_{\overline{\nabla}}u$ defines a symmetric endomorphism of $\Lambda_{-}^2(M)$.
\end{lema}
{\bf{Proof}} : \\
We change $Z$ and $U$ in $JZ$ and $JU$ respectively in (3.4) and take the sum with (3.3). Using the condition $(G_3)$, 
we get  $ \biggl [ (d_{\overline{\nabla}}u)(JZ, JU) \biggr ](X,Y)=-
\biggl [ (d_{\overline{\nabla}}u)(X,Y) \biggr ](Z,U)$. To obtain  (3.5) we change again $Z$ in $JZ$ and $U$ in $JU$ and 
use that $du_{\overline{\nabla}}(X,Y)$ is $J$-anticommuting. The rest is straighforward.
$\blacksquare$ \\ \par
This simple observation has a number of usefull consequences related to the algebraic symmetries of 
the tensor $\overline{R}$. 
\begin{coro}
Let $(M^{2n},g,J)$ be a quasi-K\"ahler manifold in the class ${\cal{AH}}_3$ and let 
$X,Y,Z,U$ be vector fields on $M$. The following holds : \\
(i) 
$$ \overline{R}(X,Y,Z,U)-\overline{R}(Z, U, X, Y)=<[\eta_X, \eta_Y]Z, U>-<[\eta_Z, \eta_U]X,Y>.$$
(ii) $\overline{R}(JX,JY),Z,U)=\overline{R}(X,Y,Z,U)$. \\
\end{coro}
{\bf{Proof}} : \\
Property (i) follows immediately from (3.2), the symmetry of Riemannian curvature operator and lemma 3.1. 
Since $\overline{R}$ is a Hermitian connection we have $\overline{R}(Z,U,JX,JY)=\overline{R}(Z,U,X,Y)$, hence 
(ii) follows from (i) and the quasi-K\"ahler condition on the tensor $\eta$.
$\blacksquare$ \\ \par
Note that the previous corollary is well known for nearly-K\"ahler manifolds (see \cite{Gray2} for instance). Also 
note that property (i) holds in fact for any almost Hermitian manifold in the class ${\cal{AH}}_3$. \par
We can have now a clearer understanding of the tensor $\tilde{R}^{\prime \prime}$ appearing in the 
Sekigawa's formula (see proposition 2.1) in terms of the torsion of the canonical Hermitian connection. 
\begin{coro}
Suppose that $(M^{2n},g,J)$ belongs to the class ${\cal{AK}}_3$. Then : 
$$ \tilde{R}^{\prime \prime}(X,Y,Z,U)=<(\overline{\nabla}_X\eta)(Y,Z)-(\overline{\nabla}_Y\eta)(X,Z),U>$$
whenever $X,Y,Z,U$ belongs to $TM$.
\end{coro}
{\bf{Proof}} :\\
This is a simple computation involving the definition of the tensor $\tilde{R}^{\prime}$ (see section 2) and lemma 
3.1. It will be left to the reader.$\blacksquare$ \\ \par
With this preliminaries in mind, we are going to show that quasi-K\"ahler manifolds in the class ${\cal{AH}}_2$ have a particularly 
nice description in terms of the torsion.
\begin{pro}
Let $(M^{2n}, g, J)$ be a quasi-K\"ahler manifold. Then $M$ satisfies condition $(G_2)$ iff
\begin{nr} \hfill
$ (\overline{\nabla}_X \eta)(Y,Z)=(\overline{\nabla}_Y \eta)(X,Z)\hfill $
\end{nr}
whenever $X,Y,Z$ are in $TM$. In particular : 
\begin{nr} \hfill
$ (\overline{\nabla}_{JX} \eta)(JY,Z)+(\overline{\nabla}_X \eta)(Y,Z)=0.\hfill $
\end{nr}
\end{pro}
{\bf{Proof}} : \\
Using (3.3) and (3.4) it is straightforward to see that condition $(G_2)$ is equivalent with 
\begin{nr} \hfill 
$ \biggl [ (d_{\overline{\nabla}}u)(Z,U) \biggr ](X,Y)=\biggl [ (d_{\overline{\nabla}}u)(JX,Y)\biggr ](JZ, U)\hfill $
\end{nr}
Expanding the right hand side of this equation and taking into account that $(M^{2n},g,J)$ is quasi-K\"ahler gives 
\begin{nr} \hfill 
$ \begin{array}{lr}
\biggl [ (d_{\overline{\nabla}}u)(Z,U) \biggr ](X,Y)=\vspace{2mm} \\
<(\overline{\nabla}_{JX}\eta)(JY,Z), U>+<(\overline{\nabla}_Y \eta)(X,Z), U>+<\eta_{T_XY}Z, U>= \vspace{2mm} \\
<(\overline{\nabla}_{JX}\eta)(JY,Z)+(\overline{\nabla}_{X}\eta)(Y,Z), U>-
\biggl [(d_{\overline{\nabla}}u)(X,Y) \biggr ](Z,U)+2< \eta_{T_XY}Z, U>.
\end{array} \hfill $
\end{nr}
But $M$ equally satisfies condition $(G_3)$, hence using (3.5) in (3.9) we obtain 
\begin{nr} \hfill 
$ \begin{array}{lr}
2\biggl [ (d_{\overline{\nabla}}u)(X,Y) \biggr ](Z,U)-2<\eta_{T_XY}Z,U>=\vspace{2mm} \\
<(\overline{\nabla}_{JX}\eta)(JY,Z), U>+<(\overline{\nabla}_X \eta)(Y,Z), U>.
\end{array} \hfill $
\end{nr}
But the right side of this equation, viewed as a tensor in $X$ and $Y$ is $J$-invariant, whilst the right 
hand side is $J$-antiinvariant. Therefore, both have to vanish and this finishes the proof
$\blacksquare$ 
\begin{rema}
From the previous proposition we deduce that quasi-K\"ahler, class ${\cal{AH}}_2$-manifolds have to satisfay the 
algebraic constraint : 
$$ <\eta_{T_XY}U,W>=<\eta_{T_UW}X,Y>.$$
Indeed, this follows immediately from (3.6) and lemma 3.1.
Note also that this relation is automatically satisfied in both nearly K\"ahler and almost K\"ahler cases. \\
\end{rema}
Using (3.1) we obtain that for a quasi-K\"ahler manifold in the class ${\cal{AH}}_2$, the difference between 
the curvature of the canonical Hermitian connection and the Riemannian curvature tensor is simply expressed by : 
\begin{nr} \hfill 
$ \overline{R}(X,Y)Z=R(X,Y)Z+[\eta_X,\eta_Y]Z-\eta_{T_XY}Z \hfill $
\end{nr}
for all vector fields $X,Y,Z$ on $M$. Note that this difference behaves as the torsion of the canonical Hermitian 
was parallel.  \par
Since $\overline{\nabla}$ is a connection with torsion it is natural to expect proposition 3.1 have further consequences. We will show 
that it is actually the case, but we have to prove first a preparatory lemma. 
\par
\begin{lema}
Let $(M^{2n},g,J)$ be a quasi-K\"ahler manifold in the class ${\cal{AH}}_2$. Then for all vector fields $X,Y,Z$ on $M$ we have : 
\begin{nr} \hfill 
$ {\sigma}_{X,Y,Z}\ \biggl ( [\overline{R}(X,Y), \eta_Z]-\eta_{\overline{R}(X,Y)Z}\biggr )+\sigma_{X,Y,Z}\ (\overline{\nabla}_{T_XY}\eta)(Z, \cdot)=0\hfill $
\end{nr}
where $\sigma$ denotes the cyclic sum. 
\end{lema}
{\bf{Proof}} : \\
Starting from $(\overline{\nabla}_Y\eta)(Z,U)=(\overline{\nabla}_Z\eta)(Y,U)$ we obtain after derivating in the direction of the vector field $X$ that : 
$(\overline{\nabla}^{2}_{X,Y} \eta)(Z,U)=(\overline{\nabla}^{2}_{X,Z} \eta)(Y,U)$. It is elementary to get then : 
$${\sigma}_{X,Y,Z}\ \biggl [ (\overline{\nabla}^{2}_{X,Y} \eta)(Z,U)-(\overline{\nabla}^{2}_{Y,X} \eta)(Z,U) \biggr ]=0. $$
Using the Ricci identity for the connection with torsion $\overline{\nabla}$ (see \cite{Besse} for instance) we arrive at : 
$${\sigma}_{X,Y,Z}\ (\overline{R}(X,Y)\eta)(Z,U)+\sigma_{X,Y,Z}\ (\overline{\nabla}_{T_XY}\eta)(Z, U)=0.$$
We end the proof by  recalling that the action of the curvature on  $\eta$ is explicitely given by : 
$$(\overline{R}(X,Y)\eta)(Z,\cdot)=[\overline{R}(X,Y), \eta_{Z}]-\eta_{\overline{R}(X,Y)Z}.$$
$\blacksquare$ \\ \par
The following result will is the starting point of the whole discussion of ${\cal{AK}}_2$ geometry in this paper.
\begin{pro}Let $(M^{2n},g,J)$ be a quasi-K\"ahler manifold in the class ${\cal{AH}}_2$. Then we have :
\begin{nr} \hfill 
$ \overline{\nabla}_{T_XY}\eta=0 \hfill $
\end{nr}
whenever  $X$ and $Y$ belong to $TM$. 
\end{pro}
{\bf{Proof}} : \\
Replacing in (3.12) $X,Y,Z$ by $JX,JY,JZ$ respectively we obtain 
$$ {\sigma}_{X,Y,Z}\ \biggl ( [\overline{R}(JX,JY), \eta_{JZ}]-\eta_{\overline{R}(JX,JY)JZ}\biggr )+\sigma_{X,Y,Z}\ (\overline{\nabla}_{T_{JX}{JY}}\eta)(JZ, \cdot)=0
$$
Now, $\eta(JU, \cdot)=-J\eta(U, \cdot)$ for all $U$ in $TM$, and it follows that 
$$ -J{\sigma}_{X,Y,Z}\ \biggl ( [\overline{R}(JX,JY), \eta_{Z}]-\eta_{\overline{R}(JX,JY)Z}\biggr )+J\sigma_{X,Y,Z}\ (\overline{\nabla}_{T_{X}{Y}}\eta)(Z, \cdot)=0.
$$
Using now corollary 3.1, (ii) we get 
$ \sigma_{X,Y,Z}\ (\overline{\nabla}_{T_{X}{Y}}\eta)(Z, \cdot)=0$, or in expanded version : 
\begin{nr} \hfill
$(\overline{\nabla}_{T_{X}{Y}}\eta)(Z, \cdot)+(\overline{\nabla}_{T_{Y}{Z}}\eta)(X, \cdot)+(\overline{\nabla}_{T_{Z}{X}}\eta)(Y, \cdot)=0. \hfill $
\end{nr}
It is clear that at $Z$ fixed, the first term of (3.14) is $J$-anti-invariant. Now, $(\overline{\nabla}_{T_{JY}{Z}}\eta)(JX, \cdot)=
-(\overline{\nabla}_{JT_{Y}{Z}}\eta)(JX, \cdot)=
(\overline{\nabla}_{T_{Y}{Z}}\eta)(X, \cdot)$ by (3.7), hence the last two terms of (3.14) are $J$-invariant. This clearly ends the proof of the proposition.
$\blacksquare$ \\ \par

\section{A first decomposition result}
In this section we will analyse some important 
geometric consequences of the proposition 3.2 established in the last section in the particular context of almost K\"ahler 
geometry. Our main object of study will be an almost K\"ahler manifold $(M^{2n},g, J)$ belonging to the 
class ${\cal{AK}}_2$. \par
Let us first set a notational convention, to be used intensively in the present and the next section
and intended to improve presentation. If $E$ and $F$ and vector subbundles of $TM$ and $Q$ is a 
tensor of type $(2,1)$, we will denote by $Q(E,F)$ (or $Q_{E}F$) the subbundle of $TM$ generated 
by elements of the form $Q(u,v)$ where $u$ belongs to $E$ and $v$ is in $F$. We will also 
denote by $<E,F>$ the product of two generic elements of $E$ and $F$ respectively. \par 
An important object associated with an almost K\"ahler manifold is its {\bf{K\"ahler nullity}}. This is the vector 
bundle $H$ over $M$ defined at a point $m$ of $M$ by $H_m=\{ v \in T_mM : \nabla_vJ=0 \}$. We also 
define ${\cal{V}}$ to be the orthogonal complement of $H$ in $TM$. Using the almost K\"ahler condition (3.1) it is easy to see that at each point of $M$ we 
have that ${\cal{V}}$ is generated by elements of the form $T_XY, X,Y$ in $TM$, in other words
$$ {\cal{V}}=T(TM, TM).$$
Hence, we have an orthogonal, $J$-invariant 
decomposition 
\begin{nr} \hfill
$TM={\cal{V}} \oplus H .\hfill $
\end{nr}
Note that, a priori, $H$ need not to have constant rank over $M$. However, this is true locally, in the following 
sense. Call a point $m$ of $M$ {\bf{regular }}
if the rank of $\eta$ attains a local maximum at $m$. Using standard continuity arguments, it follows that around each regular point, the rank of $\eta$, and 
hence that of $H$ is constant in some open subset. It is also easy to see that the set of regular point is dense 
in $M$, provided that the manifod is connected. As we are concerned with the local (in some neighbourhod of a regular 
point ) structure of ${\cal{AK}}_2$ -manifolds we can assume, without loss of generality, that $H$ has constant rank over $M$. 
This assumption will be made in the whole rest of this paper.\par
Let us examine now some elementary properties of the decomposition (4.1). 
\begin{lema}
%(i) $(\overline{\nabla}_UT)(U_1,U_2)$ belongs to $H$ for all vector fields $U, U_1,U_2$ on $TM$. \\
(i) $\overline{\nabla}_VX$ belongs to $H$ for all $V$ in ${\cal{V}}$ and $X$ in $H$.\\
(ii) $\overline{\nabla}_VW$ belongs to ${\cal{V}}$ if $V,W$ are in ${\cal{V}}$. 
\end{lema}
{\bf{Proof}} : \\
%(i) Start from the relation $<T_XY, T_ZU>=-<\eta_{T_XY}Z, U>$. Derivating in  the direction of $U_1$ we get 
%$$ \begin{array}{lr}
%<(\overline{\nabla}_{U_1}T)(X,Y), T_ZU>+<T_XY, (\overline{\nabla}_{U_1}T)(Z,U)>=\\
%-<(\overline{\nabla}_{U_1}\eta)(T_XY, Z)>-<\eta_{(\overline{\nabla}_{U_1}T)(X,Y)}Z, U>.
%\end{array}$$
%But the first term of the right hand side vanishes by (3.6) and (3.13). Moreover, by (3.5)
%$<(\overline{\nabla}_{U_1}T)(X,Y), T_ZU>=-<\eta_{(\overline{\nabla}_{U_1}T)(X,Y)}Z, U>$ hence it follows that 
%$<T_XY, (\overline{\nabla}_{U_1}T)(Z,U)>=0$ and since ${\cal{V}}$ is generated by $\{ T_XY : X,Y \in TM \}$ the 
%proof is finished. \\
(i) We know by (3.13) that $(\overline{\nabla}_V \eta)(X,U)=0$ for all $U$ in $TM$. Since $\eta_X=0$ this gives 
$\eta_{\overline{\nabla}_VX}U=0$ and the proof is finished. Now, (ii) is a straightforward consequence of (ii).
$\blacksquare$ \\ \par
As an immediate consequence we obtain our first information concerning the nature of the distributions 
${\cal{V}}$ and $H$.
%obtain using that $T(H,H)=0$ and $T({\cal{V}}, {\cal{V}}) \subseteq {\cal{V}}$ the 
%following : 
\begin{coro}
Both distributions ${\cal{V}}$ and $H$ are integrable.
\end{coro}
{\bf{Proof}} : \\
The integrability of ${\cal{V}}$ follows directly from lemma 4.1, (ii) and the fact that 
$T({\cal{V}}, {\cal{V}}) \subseteq {\cal{V}}$. To prove the integrability of $H$ we use that $(\overline{\nabla}_X\eta)(Y,U)=
(\overline{\nabla}_X\eta)(Y,U)$ for all $X,Y$ in $H$ and $U$ in $TM$ (see proposition 3.1). Since $T(H,H)=0$, this is 
equivalent to $\eta_{[X,Y]}U=0$ and this implies readily that $[X,Y]$ belongs to $H$.
$\blacksquare$ \par
Note that in the context of quasi-K\"ahler 
${\cal{AH}}_2$ manifolds it is already known \cite{Gray1} that $H$ is an integrable distribution, over each open subset 
of $M$ where it has constant rank. 
\begin{lema}
For all vector fields $V,W$ belonging to ${\cal{V}}$ and $X,Y$ in $H$ respectively we have : \\
(i) $\overline{R}(V,W)\eta=0$. \\
(ii) $[\overline{R}(X,Y), \eta_V]=\eta_{\beta_V(X,Y)}$, where $\beta_V(X,Y)=\eta_{\eta_VY}X-\eta_{\eta_VX}Y$.
\end{lema}
{\bf{Proof}} : \\
(i) We know that $\overline{\nabla}\eta$ vanishes in vertical directions (cf. 3.13). By derivation, and taking 
into account lemma 4.1, (ii) the result follows.\\ 
(ii) We use (3.12) actualized by (3.13). We have : 
$$[\overline{R}(X,Y), \eta_V]+[\overline{R}(Y,V), \eta_X] +[\overline{R}(V,X), \eta_Y]=
\eta_{\overline{R}(X,Y)V+\overline{R}(Y,V)X+\overline{R}(V,X)Y}.
$$
Now the last two terms of the first member are clearly vanishing and the use of the first Bianchi identity for $\overline{\nabla}$ yields
after a short computation to the claimed result.
$\blacksquare$ \\ \par
We are now going to obtain a first decomposition of the vector bundle ${\cal{V}}$ having good algebraic properties 
with respect to the torsion 
tensor $T$. Define a subbbundle ${\cal{V}}_0$ of ${\cal{V}}$ by setting 
$$ {\cal{V}}_0=T({\cal{V}}, {\cal{V}})$$ 
and let ${\cal{V}}_1$ be its orthogonal complement in ${\cal{V}}$. In fact, ${\cal{V}}_1$ can be considered as the K\"ahler nullity of the foliation induced 
by ${\cal{V}}$ with respect to the induced almost K\"ahler structure. More precisely, define a tensor 
$$ \hat{\eta} : {\cal{V}} \times {\cal{V}} \to {\cal{V}} \ \mbox{by} \ \hat{\eta}_vw=(\eta_vw)_{{\cal{V}}}.$$
Then we have ${\cal{V}}_1=\{ v \in {\cal{V}} : \hat{\eta}_v{\cal{V}}=0 \}$. In other words, we have 
$\eta_{{\cal{V}}_1}{\cal{V}} \subseteq H$ and this implies that $T({\cal{V}}_1, {\cal{V}}_1)=0$. Note that $\hat{\eta}$ completely determines the 
torsion over ${\cal{V}}$, that is $\hat{\eta}_vw-\hat{\eta}_wv=T(v,w)$ for all $v,w$ in ${\cal{V}}$. This follows 
from $T({\cal{V}}, {\cal{V}}) \subseteq {\cal{V}}$. \par
In the subsequent, we will assume that the subbundle ${\cal{V}}_0$ has constant rank over $M$. 
As our study is 
purely local and we have already 
assumed that $H$ has constant rank over $M$, there is no loss of generality since ${\cal{V}}_0$
has constant rank around each point of some open dense subset of $M$. 
\begin{rema}
In the subsequent we will treat the distribution ${\cal{V}}$ as it 
were an almost K\"ahler manifold with parallel 
torsion and ${\cal{V}}_1$ its K\"ahler nullity. This approach is motivated by the simple observation that the integral 
manifolds of ${\cal{V}}$ with respect to the induced metric and almost complex structure are almost K\"ahler manifolds 
with parallel torsion. Moreover, one can see that the first canonical Hermitian connection of such an integral 
manifold, coincides with the restriction of $\overline{\nabla}$.  
\end{rema}
The 
point of departure of our study is the following : 
\begin{lema}
The orthogonal decomposition ${\cal{V}}={\cal{V}}_0 \oplus 
{\cal{V}}_1$ is $J$--invariant and $\overline{\nabla}$-parallel inside ${\cal{V}}$. 
\end{lema}
{\bf{Proof}} : \\
From (3.13) we get that $\overline{\nabla}_VT=0$ for all $V$ in ${\cal{V}}$. A routine use of lemma 4.1, (ii) and of the definition of ${\cal{V_0}}$ yields now to the 
parallelism of ${\cal{V}}_0$ and hence to that of ${\cal{V}}_1$ inside ${\cal{V}}$.
$\blacksquare$ \\ \par
In the following lemma we show that the existence of such a decomposition generates strong algebraic 
restrictions involving the tensors $T$ and $\hat{\eta}$.
\begin{lema}
(i) Suppose that we have an orthogonal, $J$-invariant, decomposition ${\cal{V}}=D_1 \oplus D_2$ where the distributions $D_1$ and $D_2$ are $J$-invariant and 
$\overline{\nabla}$-parallel inside ${\cal{V}}$. Then we have : 
$$ \overline{R}(v_1,w_1,w,w_2)=-<T_{w_1}w, \hat{\eta}_{w_2}v_1>-<T_wv_1, \hat{\eta}_{w_2}v_1>$$
for all $v_1, w_1$ in $D_1$, $w$ in ${\cal{V}}$ and $w_2$ in $D_2$. \\
(ii) under the assumptions in (i) we have $\hat{\eta}_{D_2} T(D_1,D_1)=0$.\\
(iii) $\overline{R}(v,w,v_1,w_1)=0$ for all $v,w$ in ${\cal{V}}_0$ and $v_1,w_1$ in ${\cal{V}}_1$.\\
(iv) for all $v_i, 1 \le i \le 4$ in ${\cal{V}}$ we have : 
$$ \overline{R}(v_1,v_2,v_3,v_4)-\overline{R}(v_3,v_4, v_1, v_2)=
<[\hat{\eta}_{v_1}, \hat{\eta}_{v_2}]v_3, v_4>-<[\hat{\eta}_{v_3}, \hat{\eta}_{v_4}]v_1, v_2>.$$
\end{lema}
{\bf{Proof}} :  \\
We will prove (i) and (ii) in the same time. Using the first Bianchi identity for the 
connexion $\overline{\nabla}$, we get : 
$$ \begin{array}{lr}
\overline{R}(v_1, w_1, w, w_2)+\overline{R}(w_1, w, v_1, w_2)+\overline{R}(w, v_1, w_1, w_2)+\\
<T_{v_1}w_1, \eta_{w_2}w>+<T_{w_1}w, \eta_{w_2}v_1>+<T_wv_1, \eta_{w_2}w_1>=0.
\end{array}$$
Now the second and the third term below vanish since $D_1, D_2$ are $\overline{\nabla}$-parallel inside 
${\cal{V}}$. Using that 
$\overline{R}(Jv_1, Jw_1, w, w_2)=\overline{R}(v_1, w_1, w, w_2)$ and the $J$-invariance properties of the tensor 
$T$ it follows easily that $<T_{v_1}w_1, \eta_{w_2}w>=0$. Hence $<T_{v_1}w_1, \hat{\eta}_{w_2}w>=0$ and 
since $w$ in ${\cal{V}}$ was chosen arbitrary we get (ii). The proof of (i) is now straightforward.\par
To prove (iii) we take $D_1={\cal{V}}_0, D_2={\cal{V}}_1$ in (i) and use that ${\cal{V}}_1$ is the K\"ahler nullity of 
$\hat{\eta}$. Finally, for the proof of (iv) we use corollary 3.1, (i) and the fact that the tensor $\eta^H$ defined 
by $\eta^H_vw=(\eta_vw)_H$ is symmetric for all $v,w$ in ${\cal{V}}$ (this is a consequence of the fact that the 
torsion is concentrated in ${\cal{V}}$).
$\blacksquare$ \ \par
We are now able to prove our first decomposition result as follows. 
\begin{pro}
The subbundle ${\cal{V}}_0$ admits an orthogonal, $J$-invariant decomposition 
$$ {\cal{V}}_0=W_1 \oplus W_2$$
which is $\overline{\nabla}$-parallel inside ${\cal{V}}$ and has the following algebraic properties : \\
(i) $W_1=T(W_1, W_1)$ \\
(ii) $\hat{\eta}_{W_1}W_2=\hat{\eta}_{W_2}W_1=0$ \\
(iii) $\hat{\eta}_{W_2}W_2 \subseteq {\cal{V}}_1$. \\
(iv) $\hat{\eta}_{W_2}{\cal{V}}_1=W_2$.
\end{pro}
{\bf{Proof}} : \\
We define $W_1=T({\cal{V}}_0 ,{\cal{V}}_0)$ and let $W_2$ be the orthogonal complement of $W_1$ in ${\cal{V}}_0$. It is clear that $W_i, i=1,2 $ are 
$J$-invariant and $\overline{\nabla}$-parallel inside ${\cal{V}}$. Hence, we have a $\overline{\nabla}$-parallel decomposition (inside ${\cal{V}}$)
$$ {\cal{V}}=W_1 \oplus (W_2 \oplus {\cal{V}}_1). $$
Using lemma 4.4, (ii) with $D_1=W_1, D_2=W_2 \oplus {\cal{V}}_1$ we get $ \hat{\eta}_{W_1}T(W_2, {\cal{V}}_1)=0$ 
hence 
\begin{nr} \hfill 
$ \hat{\eta}_{W_1} \hat{\eta}_{W_2} {\cal{V}}_1=0 \hfill $
\end{nr}
since $\hat{\eta}_{{\cal{V}}_1}{\cal{V}}=0$. It follows that 
$$<\hat{\eta}_{W_2} {\cal{V}}_1, \hat{\eta}_{W_1}W_2>=0.$$
Now, the definition of $W_2$ implies that $<W_2, T({\cal{V}}_0, {\cal{V}}_0)>=0$ and the use of the almost 
K\"ahler condition (3.1) gives : 
\begin{nr} \hfill 
$ \hat{\eta}_{W_2} {\cal{V}}_0 \subseteq {\cal{V}}_1. \hfill $
\end{nr}
In particular, $\hat{\eta}_{W_2}W_1 \subseteq {\cal{V}}_1$. Then 
$ <\hat{\eta}_{W_2} {\cal{V}}_1, \hat{\eta}_{W_2}W_1>=0$ since the vanishing of the 
torsion on ${\cal{V}}_1$ implies that $\hat{\eta}_{W_2} {\cal{V}}_1$ is orthogonal to ${\cal{V}}_1$. We deduce that 
$<\hat{\eta}_{W_2} {\cal{V}}_1, T(W_1, W_2)>=0$, in other words $\hat{\eta}_{W_2}T(W_1,W_2)$ is orthogonal to 
${\cal{V}}_1$. But $T(W_1,W_2) \subseteq W_1$ by the definition of $W_1$ and we saw that 
$\hat{\eta}_{W_2}W_1 \subseteq {\cal{V}}_1$. It follows that 
$$\hat{\eta}_{W_2}T(W_1,W_2)=0. $$ 
Consider now the orthogonal, $J$-invariant and $\overline{\nabla}$-parallel (inside ${\cal{V}}$) decomposition ${\cal{V}}=W_2 \oplus (W_1 \oplus {\cal{V}}_1)$. Using again 
lemma 4.4, (ii) we get 
$$ \hat{\eta}_{W_2}T(W_1, W_1)=0. $$
Now, since $T(W_2, W_2)=0$ (this follows immediately from (4.3)) and by the definition of $W_1$ we have 
that $W_1$ is generated by $T(W_1,W_1)$ and $T(W_2,W_1)$ and by the previous discussion we obtain : 
$$ \hat{\eta}_{W_2} W_1=0.$$
Therefore, the second half of property (ii) is now proven .\par
Now, $\hat{\eta}_{W_2} {\cal{V}}_1$ is orthogonal to $W_1$, but we also know that it is orthogonal to 
${\cal{V}}_1$ as $T({\cal{V}}_1, {\cal{V}}_1)=0$. It follows that 
$\hat{\eta}_{W_2} {\cal{V}}_1 \subseteq W_2$. Let define now $E=\hat{\eta}_{W_2} {\cal{V}}_1$ and let 
$F$ be the orthogonal complement of $E$ in $W_2$. Then 
$\hat{\eta}_{W_2}F $ is orthogonal to ${\cal{V}}_1$ and then by (4.3) we obtain that $\hat{\eta}_{W_2}F=0$. Since 
$T(W_2, W_2)=0$ we also have 
$\hat{\eta}_{F}W_2=0$. But $\hat{\eta}_F {\cal{V}}_1 \subseteq E \subseteq W_2$ and then 
$\hat{\eta}_F {\cal{V}}_1=0$. Or $F$ is
 contained in $W_2$ hence $\hat{\eta}_FW_1=0$. We showed that 
$\hat{\eta}_F {\cal{V}}=0$ and since $F$ is contained in ${\cal{V}}_0$ it has to vanish. \par
We get that $\hat{\eta}_{W_2} {\cal{V}}_1=W_2$, proving the property (iv). Now using (4.2) we obtain that 
$\hat{\eta}_{W_1}W_2=0$ finishing the proof of (ii). Then $T(W_1, W_2)=0$ 
and this implies that $W_1=T(W_1,W_1)$. Then (i) is also proved, and (iii) is an easy consequence of (ii) and of 
vanishing of the torsion on $W_2$.
$\blacksquare$ \\ \par
We end this section with the following definition.
\begin{defi}
Let $(M^{2n},g,J)$ be in the class ${\cal{AK}}_2$. It is said to be  {\bf{special}} iff 
$H$ is of constant rank over $M$ and 
$(\nabla_{{\cal{V}}}J){\cal{V}}=H$.
\end{defi}
From an intuitive point of view, special ${\cal{AK}}_2$ manifolds are those for which the integral manifolds of the 
distribution ${\cal{V}}$, the orthogonal of the K\"ahler nullity are K\"ahler. with respect to the induced structure. 
Furthermore, the equality required in the definition forbids product with K\"ahler manifolds. 
\begin{rema} 
(i) Every $4$-dimensional almost K\"ahler manifold manifold is special in the algebraic sense of definition 4.1 on the open set where its Nijenhuis 
tensor does not vanish. This is a consequence of the fact that the vector bundle of $J$-anti-invariant $2$-forms is, in this case, of real rank $2$.\\  
(ii) If $(M^{2n},g,J)$ is a special ${\cal{AK}}_2$-manifold then it follows directly from 
the definition that $T({\cal{V}}, {\cal{V}})=0$. We also have  
$$ (\nabla_{{\cal{V}}}J)H={\cal{V}}.$$
Indeed, if $E$ is the orthogonal complement of $(\nabla_{{\cal{V}}}J)H$ in ${\cal{V}}$ then 
$\eta_{{\cal{V}}} F=0$. The vanishing of the torsion on ${\cal{V}}$ implies then
$\eta_{F} {\cal{V}}=0$. In other words, $\eta_{F}H$ is orthogonal to ${\cal{V}}$ and then 
it has to vanish. We showed that $F$ is in fact contained in the K\"ahler nullity of $(g,J)$ 
hence $F=0$ and our assertion follows.
\end{rema}
\section{Curvature properties}
In this section we will examine the curvature tensor of a local ${\cal{AK}}_2$-manifold. After proving some 
general properties we will show how the algebraic-geometric properties of $\overline{R}$ can be used to 
obtain more information about the algebraic nature of the decomposition given in proposition 4.1. Finally, this 
study will lead to the proof of theorem 1.1, which is given at the end of this section. \par
Throughout this section $(M^{2n},g,J)$ will be an almost K\"ahler manifold in the class ${\cal{AK}}_2$. All 
the notations in the previous section will be used without further comment.
\begin{lema}
Let $V_i, 1 \le i \le 3$ be in ${\cal{V}}$ and $X$ in $H$. We have : \\
(i) $\overline{R}(V_1, V_2, V_3, X)=0$ \\
(ii) $\overline{R}(X,V_1,V_2, V_3)=-<[\eta_{V_2}, \eta_{V_3}]X, V_1>$ \\
(iii) $(\overline{\nabla}_V \overline{R})(X,V_1, V_2,V_3)=0$.
\end{lema}
{\bf{Proof}} : \\
(i) follows directly from lemma 4.1, (ii) and the integrability of ${\cal{V}}$. To obtain (ii) one uses the symmetry 
property of corollary 3.1, (i). Finally, (iii) follows by derivating 
(ii) and taking into account that $\overline{\nabla}_VX$ belongs to $H$ for all $X$ in $H$ and $V$ in ${\cal{V}}$ and the fact that 
$\overline{\nabla}_V\eta=0$.
$\blacksquare$ \\ \par
We will now use the second Bianchi identity for the canonical Hermitian connection in order to get 
more information about the algebraic properties of $\nabla J$ with respect to the decomposition (4.1).
\begin{pro}
Let $X, V_i, 1 \le i \le 4$ be vector fields on $H$ and ${\cal{V}}$ respectively. We have : \\
(i) 
\begin{nr} \hfill 
$ \overline{R}(\eta_{V_2}X, V_1, V_3, V_4)-
\overline{R}(\eta_{V_1}X, V_2, V_3, V_4)=-<[\eta_{V_3}, \eta_{V_4}]X, T_{V_1}V_2>.\hfill $
\end{nr}
(ii) $(\overline{\nabla}_X \overline{R})(V_1, V_2, V_3, V_4)=0$.
\end{pro}
{\bf{Proof}} : \\
Using the second Bianchi identity we obtain 
$$ \begin{array}{lr}
(\overline{\nabla}_X \overline{R})(V_1, V_2, V_3, V_4)+(\overline{\nabla}_{V_1} \overline{R})(V_2, X,V_3, V_4)+
(\overline{\nabla}_{V_2} \overline{R})(X, V_1, V_3, V_4)+\vspace{2mm} \\
\overline{R}(T_XV_1, V_2, V_3.V_4)+\overline{R}(T_{V_1}V_2, X, V_3,V_4)+
\overline{R}(T_{V_2}X, V_1, V_3,V_4)=0.
\end{array} $$ 
Now, the second and the third terms of this equation are vanishing by lemma 5.1, (iii). It is easy to see that the first 
term is $J$-invariant in $V_1$ and $V_2$ and that all the remaining terms are $J$-anti-invariant in $V_1$ and $V_2$. 
Therefore, (ii) is proven and we obtain : $\overline{R}(T_XV_1, V_2, V_3, V_4)+\overline{R}(T_{V_1}V_2, X, V_3,V_4)+
\overline{R}(T_{V_2}X, V_1, V_3,V_4)=0.$ Since $T({\cal{V}}, {\cal{V}}) \subseteq {\cal{V}}$ it suffices now to use 
lemma 5.1, (ii) to conclude. $\blacksquare$ \\ \par
An important consequence of the equation (5.1) is : 
\begin{coro}
We have 
$$ \eta_{{\cal{V}}_1}\hat{\eta}_{{\cal{V}}_0}{\cal{V}}_0=0.$$
\end{coro} 
{\bf{Proof}} : \\
Take $V_3$ in ${\cal{V}}_0$, $V_4$ in ${\cal{V}}_1$ and $V_1, V_2$ in ${\cal{V}}$ in equation (5.1). Since ${\cal{V}}_i, i=0,1$ are 
orthogonal and $\overline{\nabla}$-parallel inside ${\cal{V}}$ we have that 
$$ <[\eta_{V_3}, \eta_{V_4}]X, T(V_1,V_2)>=0.$$
for all $X$ in $H$. Since by definition $T({\cal{V}}, {\cal{V}})={\cal{V}}_0$ it follows that $ <[\eta_{V_3}, \eta_{V_4}]X, U>=0$ for all 
$U$ in ${\cal{V}}_0$. Now, $\eta_{V_4}(\eta_{V_3}X)$ is in $H$ since $V_4$ is in ${\cal{V}}_1$ and 
$\eta_{V_3}X$ in ${\cal{V}}$ hence $ <\eta_{V_4}X, \hat{\eta}_{V_3}U>=-<\eta_{V_3}\eta_{V_4}X, U>=0$. In our 
notations, $<\eta_{{\cal{V}}_1}H, \hat{\eta}_{{\cal{V}}_0}{\cal{V}}_0>=0$. Then we get that 
$\eta_{{\cal{V}}_1}\hat{\eta}_{{\cal{V}}_0}{\cal{V}}_0$ is orthogonal to $H$. Or the definition of ${\cal{V}}_1$ ensures that 
$\eta_{{\cal{V}}_1}{\cal{V}}$ is contained in $H$ and our result follows.
$\blacksquare$ \\ \par
Using the previous corollary and equation (5.1) as main tools, we will proceed now to the refinement of the decomposition given in proposition 
4.1. More precisely, our immediate objective will be to show that the space $W_2$ occuring in the decomposition 
given in proposition 4.1, must vanish. 
To  proceed, consider the decomposition 
$$ {\cal{V}}_0=W_1 \oplus W_2 $$
and define $E^{\prime}=\hat{\eta}_{W_2}W_2 \subseteq {\cal{V}}_1$ and let $E$ be the orthogonal complement of $E^{\prime}$ in ${\cal{V}}_1$. Obviously, the 
decomposition ${\cal{V}}_1=E \oplus E^{\prime}$ is $J$-invariant and $\overline{\nabla}$-parallel inside 
${\cal{V}}$. \par
Using corollary 5.1 we obtain some preliminary algebraic information as follows. 
\begin{lema} We have : 
$$ \eta_{E^{\prime}}H \subseteq W_2. $$
\end{lema}
{\bf{Proof}} : \\
By corollary 5.1, the fact that $T(W_1, W_1)=W_1$ and the definition of $E^{\prime}$ we obtain easily that 
\begin{nr} \hfill
$ \eta_{{\cal{V}}_1}W_1=0 \hfill $ 
\end{nr}
and 
\begin{nr} \hfill
$ \eta_{{\cal{V}}_1} E^{\prime}=0. \hfill $
\end{nr}
The second equation gives us $\eta_E E^{\prime}=\eta_{E^{\prime}}E^{\prime}=0$. Then we have that 
$\eta_{E^{\prime}} E$ is contained in $T(E, E^{\prime})=0$ (since $T({\cal{V}}_1, {\cal{V}}_1)=0$) and thus it vanishes, 
showing that $\eta_{E^{\prime}} {\cal{V}}_1=0$. Hence $\eta_{E^{\prime}}H$ is orthogonal to ${\cal{V}}_1$ and by 
the first equation it is also orthogonal to $W_1$ and our result follows.
$\blacksquare$ \\ \par
We will need now one more auxiliary lemma. 
\begin{lema}
We have : \\
(i) 
$$ \overline{R}(v_1,w_1,v_0, w_0)=-<[\hat{\eta}_{v_0}, \hat{\eta}_{w_0}]v_1, w_1>$$
for all $v_1,w_1$ in ${\cal{V}}_1$ and $v_0,w_0$ in ${\cal{V}}_0$. \\
(ii) Define a symmetric tensor $\hat{r} : {\cal{V}}_0 \to {\cal{V}}_0$ by setting 
$$ <\hat{r}v_0,w_0>= \sum \limits_{v_k \in {\cal{V}}_1}^{}<\hat{\eta}_{v_0}v_k, \hat{\eta}_{w_0}v_k>$$
for all $v_0, w_0$ in ${\cal{V}}_0$ (here $\{v_k \}$ is an arbitrary local orthonormal 
basis in ${\cal{V}}_1$). Then $\hat{r}$ preserves $W_2$ and the restriction 
of $\hat{r}$ to $W_2$ has no kernel.\\
(iii) Let $U$ be in ${\cal{V}}$ such that 
$$ \overline{R}(U,V_1,V_2,V_3)=0$$
for all $V_1$ in $E^{\prime}$ and $V_2, V_3$ in ${\cal{V}}_1$. Then $U \perp E^{\prime}$.
\end{lema}
{\bf{Proof}} : \\
(i) It suffices to apply lemma 4.4, (iii) and (iv). \\
(ii) That $W_2$ is preserved by $\hat{r}$ follows directly by the fact that $W_2$ is $\overline{\nabla}$-parallel 
inside ${\cal{V}}$. 
Let $v$ in $W_2$ be such $\hat{r}V=0$. Then the definition of $\hat{r}$ implies 
directly that $\hat{\eta}_v {\cal{V}}_1=0$. It follows that $\hat{\eta}_vW_2$ is orthogonal to ${\cal{V}}_1$ and we know 
that it is also contained in ${\cal{V}}_1$, since $v$ belongs to $W_2$ (cf. proposition 4.1, (iii)). Thus 
$\hat{\eta}_v W_2=0$, and 
again the fact that $v$ belongs to $W_2$ yields $\hat{\eta}_v W_1=0$ (see proposition 4.1, 
(ii)). Hence $\hat{\eta}_v {\cal{V}}=0$, $v$ is in 
$W_2 \subseteq {\cal{V}}_0$ and this clearly implies the vanishing of $v$. \\
(iii) Using the symmetry formula of lemma 4.4, 
(iv) we obtain that 
$$\overline{R}(V_2,V_3, V_1,U)=0$$ 
for all $V_2,V_3$ in ${\cal{V}}_1$ and $V_1$ in $E^{\prime}$. 
%Let $U_1$ be the orthogonal projection of $U$ on $E^{\prime}$. As $E^{\prime}$ is parallel inside ${\cal{V}}$ we obtain 
%that 
%$$ \overline{R}(V_1,V_2, U_1, V_1)=0$$
%for all $V_1,V_2$ in ${\cal{V}}_1$ and $V_1$ in $E^{\prime}$. 
%Consider now an element $V_1$ of $E^{\prime}$, of the 
%for $\hat{\eta}_{v_0}w_0$ with $v_0, w_0$ in $W_2$. 
Let now $v_0,w_0$ be in $W_2$.  If $\{ v_k\}$ is an orthonomal basis in ${\cal{V}}_1$ then 
by lemma 4.2, (i) we have that : 
\begin{nr} \hfill
$ \overline{R}(v_k,Jv_k) \eta_{v_0}w_0=\eta_{\overline{R}(v_k,Jv_k)v_0}w_0+
\eta_{v_0}(\overline{R}(v_k,Jv_k)w_0). \hfill $
\end{nr}
Now note that using point (i), one easily finds that $\sum \limits_{v_k \in {\cal{V}}_1}
\overline{R}(v_k,Jv_k)v=-2(\hat{r}J)v$ for all $v$ in ${\cal{V}}_0$. With this in mind, we project (5.4) on 
${\cal{V}}$ and sum over $k$ to find : 
$$ \sum \limits_{v_k \in {\cal{V}}_1}^{}\overline{R}(v_k,Jv_k) \hat{\eta}_{v_0}w_0=2J \biggl [\hat{\eta}_{\hat{r}v_0}w_0+
\hat{\eta}_{v_0}(\hat{r}w_0)].$$
But elements of the form $\hat{\eta}_{v_0}w_0$ generates, by definition, $E^{\prime}$ hence 
taking  the scalar product with $U$ yields to $<\hat{\eta}_{\hat{r}v_0}w_0+
\hat{\eta}_{v_0}(\hat{r}w_0), JU>=0$ for all $v_0, w_0$ elements of $W_2$. Since $\hat{r}$ is positive 
without kernel on $W_2$, by eventually considering its spectral decomposition 
we deduce that $JU$ (and then $U$)  is orthogonal to $\hat{\eta}_{W_2}W_2=E^{\prime}$ and the proof is finished.
$\blacksquare$ \\ \par
Based on these preparations we are now able to show the following.
\begin{pro}
The subbundle ${\cal{V}}$ admits a orthogonal, 
$J$-invariant and $\overline{\nabla}$-parallel (inside ${\cal{V}}$)-decomposition : 
$$ {\cal{V}}={\cal{V}}_0 \oplus {\cal{V}}_1$$
with ${\cal{V}}_0=T({\cal{V}}_0, {\cal{V}}_0)$ and $\eta_{{\cal{V}}_1}{\cal{V}}_0=0$.
\end{pro}
{\bf{Proof}} : \\
We are going to show first that we must have $E^{\prime}=0$. Let us consider $v_0$ in $W_2$, $w_0$ in 
$E^{\prime}$ and $v_1,v_2$ in ${\cal{V}}_1$ as well as $X$ in $H$. Using proposition 5.1, (i) we obtain 
$$ \overline{R}(\eta_{v_0}X, w_0, v_1,v_2)-
\overline{R}(\eta_{w_0}X,v_0, v_1,v_2)=-<[\eta_{v_1}, \eta_{v_2}] X, T_{w_0}v_0>.$$
But $\eta_{w_0}X$ belongs (see lemma 5.2) to $W_2 \subseteq {\cal{V}}_0$, hence the second term 
of the left hand side vanishes by lemma 4.4, (iii), as well as the right hand side, since $T_{w_0}v_0$ is in ${\cal{V}}$ whilst 
$[\eta_{v_1}, \eta_{v_2}] X$ is in $H$. We found that : 
$$ \overline{R}(\eta_{v_0}X, w_0, v_1,v_2)=0$$
for all $w_0$ in $E^{\prime}$ and all $v_1,v_2$ in ${\cal{V}}_1$. Applying lemma 5.3, (iii), we obtain that 
$\eta_{v_0}H$ is orthogonal to $E^{\prime}$, that is in our notations $<\eta_{W_2}H, E^{\prime}>=0$. It follows
that $\eta_{W_2}E^{\prime} \subseteq {\cal{V}}=W_1 \oplus W_2 \oplus {\cal{V}}_1$. But 
$\hat{\eta}_{W_2}W_1=0, \hat{\eta}_{W_2}{\cal{V}}_1=W_2$  (cf. proposition 4.1, (ii) and 
(iv)) hence $\eta_{W_2}E^{\prime} \subseteq W_2$. Or $T(E^{\prime}, W_2) \subseteq {\cal{V}}$ hence it follows 
that $\eta_{E^{\prime}}W_2 \subseteq {\cal{V}}$. Since 
$E^{\prime}$ lies in ${\cal{V}}_1$, the definition of ${\cal{V}}_1$ implies $\eta_{E^{\prime}}W_2 \subseteq H$ and 
we obtain that $\eta_{E^{\prime}}W_2=0$.
%and this implies $\eta_{E^{\prime}}W_2 \subseteq W_2$. As 
%$T(E^{\prime}, W_2) \subseteq {\cal{V}}$ we find that $\eta_{E^{\prime}}W_2=0$. 
By means of lemma 5.2 this also implies $\eta_{E^{\prime}}H=0$. But using (5.2) and (5.3) we obtain (because $E^{\prime}$ is contained in 
${\cal{V}}_1$) that $\eta_{E^{\prime}}W_1=0$ and $\eta_{E^{\prime}}{\cal{V}}_1=0$. This means that 
$E^{\prime}$ is contained in $H$, the K\"ahler nullity of $(g,J)$ and we obtain that $E^{\prime}=0$. \par
Now, $\hat{\eta}_{W_2}W_2$ must vanish and it follows that $\hat{\eta}_{W_2}{\cal{V}}_1$ is orthogonal to 
$W_2$. But we already know (see proposition 4.1, (iv)) that the last mentioned spaces are in fact equal, and then 
$W_2=0$. It follows that $W_1={\cal{V}}_0$ and the proof is now finished, where the last assertion follows by (5.2).
$\blacksquare$ \\ \par
We will investigate now the geometrical properties of the decomposition in the proposition 5.2, properties that will 
lead, once again, to a better understanding of its algebraic structure. \par
Let we introduce the {\it{configuration}} tensor 
$A : H \times H \to {\cal{V}}$ by setting : 
$$ \overline{\nabla}_XY=\tilde{\nabla}_XY+A_XY$$ 
for all $X,Y$ in $H$, where $\tilde{\nabla}_XY$ denotes the orthogonal projection of $\overline{\nabla}_{X}Y$ 
on $H$. The tensor $A$ is precisely the obstruction to the distribution $H$ to be $\overline{\nabla}$-
parallel. In a similar way, we define $B : H \times {\cal{V}} \to H$ by 
$$ \overline{\nabla}_XV=\tilde{\nabla}_XV+B_XV.$$
Because the connection $\overline{\nabla}$ is metric we have $<B_XV, Y>=-<V,A_XY>$ for all $X,Y$ in $H$ and 
$V$ in ${\cal{V}}$. Since $H$ is integrable we have that $A$ is symmetric, that is $A_XY=A_YX$. It is also obvious 
that $ [A_X, J]=0$ for all $X$ in $H$. 
\begin{lema}
(i) The distribution ${\cal{V}}_0$ is $\overline{\nabla}$-parallel.\\
(ii) Moreover we have that $\eta_{{\cal{V}}_0}{\cal{V}}_1=0$ and $\eta_{{\cal{V}}_1}H={\cal{V}}_1$. 
\end{lema}
{\bf{Proof}} : \\
(i) We begin by showing that the operator $B_X, X$ in $H$ must vanish on ${\cal{V}}_0$. Indeed, let us recall that 
$(\overline{\nabla}_X \eta)(V, \cdot)=(\overline{\nabla}_V \eta)(X, \cdot)=0$ for all $X$ in $H$ and $V$ in ${\cal{V}}$ 
(see (3.6) and (3.13)). If $W$ is in ${\cal{V}}$ it follows then easily that 
$(\overline{\nabla}_XT)(V,W)=0$. In other words 
$$ \tilde{\nabla}_X(T_VW)+B_X(T_VW)=T(\overline{\nabla}_XV,W)+T(V,\overline{\nabla}_XW)$$
belongs to ${\cal{V}}$ and this implies the vanishing of $B_X$ on ${\cal{V}}_0$. Consider now $V_0, W_0$ in 
${\cal{V}}_0$. We have, for any $X$ in $H$ : 
$$ \begin{array}{rr}
\overline{\nabla}_X (T_{V_0}W_0)=T(\overline{\nabla}_XV_0, W_0)+T(V_0, \overline{\nabla}_XW_0)\\
=T(\tilde{\nabla}_XV_0, W_0)+T(V_0, \tilde{\nabla}_XW_0).
\end{array}$$
But ${\cal{V}}_0=T({\cal{V}}_0, {\cal{V}}_0)=T({\cal{V}}_0, {\cal{V}})$ and ${\cal{V}}_0$ is $\overline{\nabla}$-parallel 
insider ${\cal{V}}$, hence (i) is proven. \\
(ii) We will show first that $\eta_{{\cal{V}}_1}H={\cal{V}}_1$. Since $\eta_{{\cal{V}}_1}{\cal{V}}_0=0$ we have that 
$\eta_{{\cal{V}}_1}H \subseteq {\cal{V}}_1$. Consider the decomposition ${\cal{V}}_1=E \oplus F$ with 
$\eta_{{\cal{V}}_1}H=E$ and $F$ the orthogonal complement of $E$ in ${\cal{V}}_1$. From the definition of 
$F$ it follows that $\eta_{{\cal{V}}_1}F$ is orthogonal to $H$ and hence it vanishes (recall 
that $\eta_{{\cal{V}}_1} {\cal{V}}_1 \subseteq H$). Since $T({\cal{V}}_1, {\cal{V}}_1)=0$ it also follows that 
$\eta_F {\cal{V}}_1=0$. This implies that $\eta_F H$, which a subspace of ${\cal{V}}_1$, is orthogonal to 
${\cal{V}}_1$, and hence $\eta_FH=0$. Finally since $\eta_F {\cal{V}}_0=0$ ($F$ lies in ${\cal{V}}_1$) we get 
that $F$ is contained in the K\"ahler nullity of $(g,J)$ and then, of course, $F$=0. \par
Now, by (i) we deduce that ${\cal{V}}_1 \oplus H$ is a $\overline{\nabla}$-parallel distribution. Using an argument 
similar to that of lemma 4.4, (ii) for the $\overline{\nabla}$-parallel decomposition 
$TM={\cal{V}}_0 \oplus ({\cal{V}}_1 \oplus H)$ we find that 
$$ \eta_{{\cal{V}}_0} T({\cal{V}}_1, H)= \eta_{{\cal{V}}_0} \eta_{{\cal{V}}_1}H=0.$$ 
Combining this with $\eta_{{\cal{V}}_1}H={\cal{V}}_1$ finishes the proof of the lemma.
$\blacksquare$ \\ \par
The last step before proving the splitting theorem 1.1, consists in investigating reducibility properties of 
the K\"ahler nullity $H$. We need to introduce some notations. For every $X$ in $H$ define a linear map : 
$$ \gamma_X : {\cal{V}}_1 \to {\cal{V}}_1 \ \mbox{by} \ \gamma_XV=\eta_VX .$$
The maps $\gamma_X$ are in relation with the curvature of $H$ (with respect to the canonical Hermitian connection), 
as showed in the following lemma.
\begin{lema} Let $X,Y,Z$ be in $H$ and $V,W$ in ${\cal{V}}_1$. We have : 
$$ \overline{R}(X,Y, \eta_VW, Z)=\overline{R}(\gamma_ZV,W,X,Y)+<[[\gamma_X, \gamma_Y],\gamma_Z]V,W>$$
\end{lema}
{\bf{Proof}} : \\
From lemma 4.2, (ii), we know that $[\overline{R}(X,Y), \eta_V]W=\eta_{\beta_V(X,Y)}W$, where, in our present 
notations $\beta_V(X,Y)=\eta_{\gamma_YV}X-\eta_{\gamma_XV}Y=[\gamma_X, \gamma_Y]V$. Taking the scalar 
product with $Z$ yields after an easy calculation to the wanted result.
$\blacksquare$ \\ \par
\begin{lema}
Suppose that $H_1$ is a $J$-invariant, $\overline{\nabla}$-parallel distribution contained in $H$. Let 
$H_2$ be the orthogonal complement of $H_1$ in $H$. Then the spaces $\eta_{{\cal{V}}_1}H_1$ and 
$\eta_{{\cal{V}}_1}H_2$ are mutually orthogonal.
\end{lema}
{\bf{Proof}} : \\
Let $X_1$ and $X_2$ be in $H_1$ and $H_2$ respectively. Then the parallelism of $H_1$, together with the symmetry 
property of $\overline{R}$ (see corollary 3.1, (i)) ensures that 
$\overline{R}(X_1,X_2, \eta_VW, Z)=\overline{R}(\gamma_ZV, W, X_1,X_2)=0$ for all $V,W$ in ${\cal{V}}_1$ and 
$Z$ in $H$. Then, by the previous lemma we obtain 
$$ [[\gamma_{X_1}, \gamma_{X_2}],\gamma_Z]=0$$
for all $Z$ in $H$. Taking $Z=X_1$ we find that 
$$ \gamma_{X_1}^2 \gamma_{X_2}+\gamma_{X_2}\gamma_{X_1}^2=
2 \gamma_{X_1} \gamma_{X_2} \gamma_{X_1}.$$
We change now $X_2$ in $JX_2$ in the previous equation and take into account that 
$\gamma_{JX}=\gamma_XJ=-J\gamma_X$. It follows that 
$$ \gamma_{X_1}^2 \gamma_{X_2}+\gamma_{X_2}\gamma_{X_1}^2=
-2 \gamma_{X_1} \gamma_{X_2} \gamma_{X_1}$$
hence we must have 
$$\gamma_{X_1}^2 \gamma_{X_2}+\gamma_{X_2}\gamma_{X_1}^2=
 \gamma_{X_1} \gamma_{X_2} \gamma_{X_1}=0.$$
This easily implies that $\gamma_{X_1}^3 \gamma_{X_2}=0$ and since $\gamma_{X}$ is a symmetric operator 
for all $X$ in $H$ (as a consequence of $T({\cal{V}}, {\cal{V}}) \subseteq {\cal{V}}$) 
we get that $\gamma_{X_1} \gamma_{X_2}=0$. But this fact is equivalent to the orthogonality 
of the spaces $\eta_{{\cal{V}}_1}H_1$ and 
$\eta_{{\cal{V}}_1}H_2$, and the proof is now finished.
$\blacksquare$ \\ 
$\\$ 
{\bf{Proof of theorem 1.1}} : \\
Let us define the distribution $H_1$ to be $H_1=(\eta_{{\cal{V}}_0}{\cal{V}}_0)_H$ where the subscript denotes 
orthogonal projection. We have then $\eta_{{\cal{V}}_0}{\cal{V}}_0={\cal{V}}_0 \oplus H_1$. As ${\cal{V}}_0$ is 
$\overline{\nabla}$-parallel, so is $\eta_{{\cal{V}}_0}{\cal{V}}_0$, hence $H_1$ must be $\overline{\nabla}$-parallel. Define 
now ${\cal{W}}_i=\eta_{{\cal{V}}_1}H_i, 1 \le i \le 2$. Then using lemma 5.6 and the fact that $\eta_{{\cal{V}}_1} H=
{\cal{V}}_1$ we obtain a $J$-invariant and orthogonal 
decomposition 
$${\cal{V}}_1={\cal{W}}_1 \oplus {\cal{W}}_2.$$
The orthogonality of ${\cal{W}}_1$ and ${\cal{W}}_2$ ensures, in the standard way, that $\eta_{{\cal{W}}_1}{\cal{W}}_2=0$ and that 
$\eta_{{\cal{W}}_1}H_2=\eta_{{\cal{W}}_2}H_2=0$. Let us show now that ${\cal{W}}_1$ is a $\overline{\nabla}$-parallel distribution. 
Let $X$ be in $H$ and $U,X_1$ be in ${\cal{W}}_1$ and $H_1$ respectively. Then : 
$$\overline{\nabla}_X (\eta_UX_1)=\eta_U \overline{\nabla}_XX_1+\eta_{\overline{\nabla}_XU}X_1$$
belongs to $\eta_{{\cal{W}}_1}H_1+\eta_{TM}X_1={\cal{W}}_1$, were we used the $\overline{\nabla}$-parallelism of 
$H_1$. In the same way it can be showed that  ${\cal{W}}_1$ is $\overline{\nabla}$-parallel inside ${\cal{V}}$. We get 
a $\overline{\nabla}$-parallel decomposition : 
$$ TM=({\cal{V}}_0 \oplus {\cal{W}}_1 \oplus H_1) \oplus ({\cal{W}}_2 \oplus H_2).$$
Using the behaviour of the tensor $\eta$ with respect to this decomposition proven previously 
it follows that this splitting 
is in fact $\nabla$-parallel. From the discussion below it is now clear that the first factor gives rise 
to an ${\cal{AK}}$-manifold with parallel torsion and the second to a special ${\cal{AK}}_2$-manifold.  
$\blacksquare$ \\ \par

\section{On parallel torsion}
In this section we will consider an almost K\"ahler manifold $(M^{2n},g,J)$ whose first canonical connection has parallel torsion. Then proposition 3.1 tells us that 
$(M^{2n},g,J)$ naturally belongs to the class ${\cal{AK}}_2$. Our aim here is to show that in this setting that metric cannot be Einstein if the 
$J$ is not integrable. Along the way we will also obtain some information about the holonomy of the canonical Hermitian connection. We will 
need first some preliminaries. At first, we define a symmetric tensor $r : TM \to TM$ by $<rX,Y>=\sum \limits_{i=1}^{2n}<(\nabla_{e_i}J)X, (\nabla_{e_i}J)Y>$ for all 
$X$ and $Y$ in $TM$, where 
$\{ e_i, 1 \le i \le 2n \}$ is an arbitrary local orthonomal basis of $TM$. 
\begin{lema} For any ${\cal{AK}}$-manifold with parallel torsion we have : \\
(i) $$\nabla^{\star} \nabla \omega=\Phi^0$$ 
where the $2$-form $\Psi$ is defined by $\Psi(X,Y)=<(rJ)X,Y>$. \\
(ii) If $\{ e_i, 1 \le i \le 2n\}$ is an arbitrary local orthonormal basis in $TM$ we define the $2$-form $\overline{\rho}$ by  setting 
$ \overline{\rho}=\sum \limits_{i=1}^{2n}<\overline{R}(e_i,Je_i) \cdot, \cdot >$. Then we have : 
$$ \overline{\rho}=2\rho+\frac{1}{2} \Psi.$$
\end{lema}
{\bf{Proof}} : \\
(i) Let us recall the fact that $(\nabla_X \omega)(Y,Z)=<(\nabla_XJ)Y,Z>$ for all $X,YZ$ in $TM$. Using this and the parallelism of the torsion, a simple algebraic computation which 
will left to the reader yields to the desired result.\\
(ii) Follows immediately from formula (3.11) and (2.1).
$\blacksquare$ \\ \par
It is an easy exercise to show that the Ricci tensor is $J$-invariant in the presence of parallel torsion. Furthermore, from lemma 6.1, (i) and (2.1) it also follows that the 
Ricci-$\star$ tensor is $J$-invariant. \par 
Now, our key ingredient for proving theorem 1.2 consists in the following lemma, whose first part is Sekigawa's formula in the ${\cal{AK}}_2$ case and 
whose second part is a relation complementary to Sekigawa's coming from the parallelism of the torsion. 
\begin{lema} We have : \\
(i) 
$$ 4<\rho, \Phi-\Psi>=\vert \Phi \vert^2+\vert \Psi \vert^2$$
(ii) 
$$ 4<\rho, \Phi+\Psi>+<\Psi, \Phi+\Psi>=0.$$
\end{lema}
{\bf{Proof}} : \\
(i) This is an immediate consequence of Sekigawa's formula (see proposition 2.1), actualized in the parallel torsion case. Indeed, we must have 
$\tilde{R}^{\prime \prime}=0$ (see corollary 3.2). It also clear that $s^{\star}-s$ is a constant function. Then, using lemma 6.1, (i) in Sekigawa's formula we get the desired result.\\
(ii) Since the torsion is parallel, we have that $\overline{R}(X,Y).\eta=0$ for all $X,Y$ in $TM$. It follows that $\overline{\rho}.\eta=0$ and taking the scalar product with 
$J\eta$ we obtain after an easy computation that $<\overline{\rho}, \Phi+\Psi>=0$. It suffices now to use lemma 6.1, (ii) to conclude. 
$\blacksquare$ \\ \par
We can prove the following proposition, which is nothing else that theorem 1.2 in the introduction.
\begin{pro}
Let $(M^{2n},g,J)$ be almost K\"ahle with parallel torsion. Then : \\
(i) If $g$ is Einstein then $J$ is integrable ; \\
(ii) If $J$ is not integrable the connection $\overline{\nabla}$ has real reducible holonomy.
\end{pro}
{\bf{Proof}} : \\
We prove both assertions in the same time. Let us suppose that we have 
\begin{nr} \hfill
$2<\rho, \Phi>=<\rho, \Psi> \hfill $ 
\end{nr}
and prove that $J$ is integrable. Using (7.1), the relations in lemma 7.1 become : 
$$\begin{array}{lr}
-2<\rho, \Psi>=\vert \Phi \vert^2+\vert \Psi \vert^2 \\
6<\rho, \Psi>+\vert \Psi \vert^2+<\Psi, \Phi>=0.
\end{array} $$
We deduce that $3\vert \Phi \vert^2+2\vert \Psi \vert^2=<\Phi, \Psi>$. Since $<\Phi, \Psi>\le 
\vert \Phi \vert \cdot \Psi \vert$ we have clearly that $\Psi=\Phi=0$, that is $(g,J)$ is a K\"ahler structure. \par
Now, if the manifold is Einstein, (6.1) is clearly satisfied, hence (i) is proven. To prove (ii), suppose that 
$\overline{\nabla}$ has irreducible holonomy. Then the $\overline{\nabla}$-parallel forms $\Phi, \Psi$ must 
be multiples of $\omega$ hence $2\Phi=\Psi$ so (6.1) is again satisfied.
$\blacksquare$ \\ \par
In the same vein, one can also have integrability results in terms of the Hermitian Ricci tensor $\overline{\rho}$. 
\begin{pro}Let $(M^{2n},g,J)$ be almost K\"ahler with parallel torsion. If there exists a real constant $\lambda$ such 
that $\overline{\rho}=\lambda J$ then $(g,J)$ is a K\"ahler structure.
\end{pro}
The proof will be ommited since completely analogue to the of proposition 6.1. \\ 
$\\$
{\bf{Acknowledgements}} : The author would like to thank T. Draghici for many interesting discussions on 
almost K\"ahler geometry and also for pointing out some imprecisions in an earlier version of this paper.

\begin{flushright}
Paul-Andi Nagy \\
Institut de Math\'ematiques \\
rue E. Argand 11, CH-2007, Neuch\^atel \\ 
email : Paul.Nagy@unine.ch
\end{flushright}

\normalsize
\begin{thebibliography}{99}
\bibitem{Ivanov}
ALEXANDROV. B, GRANTCHAROV. G, IVANOV. S, \textit{Curvature properties of twistor spaces of quaternionic K\"ahler manifolds}, 
J. Geom. {\bf{69}} (1998), 1-12.
\bibitem{Apo3}
V. APOSTOLOV, T. DRAGHICI, KOTSCHICK, D. \textit{An integrability theorem for almost K\"ahler $4$-manifolds}, C.R.Acad.Sci. Paris{\bf{329}}, 
s\'er. I (1999), 413-418.
\bibitem{Apo}
APOSTOLOV, V., CALDERBANK, D., GAUDUCHON, P., \textit{The geometry of weakly 
selfdual K\"ahler surface }, to appear in Compositio Math. 
\bibitem{Apo2}
V. APOSTOLOV, T. DRAGHICI, A. MOROIANU, \textit{A splitting theorem for K\"ahler manifolds with constant eigenvalues 
of the Ricci tensor}, Int. J. Math. {\bf{12}} (2001), 769-789.
\bibitem{Apo1}
V. APOSTOLOV, J. ARMSTRONG, T. DRAGHICI 
\textit{Local rigidity of certain classes of almost K\"ahler $4$-manifolds}, Ann. Glob. Anal. Geom. {\bf{21}} (2002), 151-176.
\bibitem{Apo4}
V. APOSTOLOV, J. ARMSTRONG, T. DRAGHICI, \textit{Locals models and integrability of certain almost K\"ahler $4$-manifolds}, Math. Ann., 
to appear.
\bibitem{arm1}
ARMSTRONG, J. \textit{ On four-dimensional almost K\"ahler manifolds}, Quart. J. Math. Oxford Ser. (2) {\bf{48}} (1997), 405-415.
\bibitem{arm2}
ARMSTRONG, J. \textit{An Ansatz for Almost-K\"ahler, Einstein $4$-manifolds }, J. reine angew. Math. {\bf{542}} (2002), 53-84.
\bibitem{Besse}
A. L. BESSE, \textit{Einstein manifolds}, Springer Verlag, 1986.
\bibitem{Davidov}
J. DAVIDOV, O. MUSKAROV, \textit{Twistor spaces with Hermitian Ricci tensor}, Proc. Amer. Math. Soc. {\bf{109}} (1990), no.4, 1115-1120.
%\bibitem{Kowalski}
%O. KOWALSKI, \textit{Generalized Symmetric spaces}, LNM {\bf{805}}, 1980.
\bibitem{Golberg}
S. I. GOLDBERG, \textit{Integrability of almost K\"ahler manifolds}, Proc. Amer. Math. Soc. {\bf{21}} (1969), 96-100.
\bibitem{Gray1}
A. GRAY, \textit{Curvature identities for Hermitian and almost Hermitian manifolds}, T\^ohoku Math. J. 
{\bf{28}} (1976), 601-612.
\bibitem{Gray2}
A. GRAY, \textit{The structure of Nearly K\"ahler manifolds}, Math. Ann. {\bf{223}} (1976), 233-248.
\bibitem{Nagy2}
NAGY, P-A, \textit{Nearly K\"ahler geometry and Riemannian foliations }, Asian J. Math., {\bf{6}} no. 3 (2002), 481-504, to appear. 
\bibitem{Nur}
NUROWSKI, P., PRZANOWSKI, M, \textit{A four-dimensional example of Ricci flat metric 
admitting K\"ahler non-K\"ahler structure }, Classical Quantum Gravity {\bf{16}} (1999), no.3, L9-L13.
\bibitem{Seki2}
OGURO. T, SEKIGAWA, K. \textit{Four dimensional almost K\"ahler Einstein and $\star$-
Einstein manifolds}, Geom. Dedicata {\bf{69}} (1998), no.1, 91-112.
%\bibitem{Tri2}
%F. TRICERRI, L. VANHECKE, \textit{Curvature tensors on almost Hermitian manifolds}, Trans. Amer. 
%math. Soc. {\bf{267}} (1981), 365-398.
%\bibitem{Tri}
%F. TRICERRI, L. VANHECKE, \textit{On a theorem of Kiricenko relating to $3$-symmetric spaces}, Riv. Mat. Uni. Parma
\bibitem{Olszak}
Z. OLSZAK, \textit{A note on almost K\"ahler manifolds}, Bull. Acad. Polon. Sci. {\bf{XXVI}} (1978), 199-206.
\bibitem{Seki1}
SEKIGAWA, K, \textit{on some compact Einstein almost K\"ahler manifolds}, J. Math. Soc. Japan {\bf{36}} 
(1987), 677-684.
\end{thebibliography}
\end{document}